\documentclass[12pt]{article} 
\usepackage{hyperref}
\usepackage{amssymb}
\usepackage{amsfonts}
\usepackage{amsmath}

\usepackage{enumerate}

\newtheorem{theorem}{Theorem}

\newtheorem{proposition}[theorem]{Proposition}
\newtheorem{lemma}[theorem]{Lemma}
\newtheorem{corollary}[theorem]{Corollary}
\newtheorem{remark}[theorem]{Remark}

\newcommand{\proof}{ \noindent{\it Proof:\ \ }}
\def\qed{\ifhmode\unskip\nobreak\fi\ifmmode\ifinner\else\hskip5 pt
\fi\fi\hbox{\hskip5 pt \vrule width4 pt height6 pt depth1.5 pt
\hskip 1pt }}
\newcommand{\po}{{\hspace*{-1ex}}{\bf .  }}

\let\oldmarginpar\marginpar
\renewcommand\marginpar[1]{${}^\clubsuit$\oldmarginpar[\raggedleft\scriptsize\sf #1]{\raggedright\scriptsize\sf #1}}
\usepackage{titlesec}
\titleformat{\section}{\normalfont\sffamily\Large}{\thesection.}{0.7em}{}
\titleformat{\subsection}{\normalfont\sffamily\filcenter}{}{0.7em}{}

\newcommand{\immr}[2]{$f: M^{#1} \rightarrow \R^{#1 + #2}$}
\newcommand{\immrw}[2]{$\hat{f}: M^{#1} \rightarrow \R^{#1 + #2}$}

\newcommand{\cref}[1]{Corollary~\ref{#1}}
\newcommand{\csref}[1]{Corollaries~\ref{#1}}

\newcommand{\lref}[1]{Lemma~\ref{#1}}
\newcommand{\pref}[1]{Proposition~\ref{#1}}

\newcommand{\sref}[1]{Section~\ref{#1}}
\newcommand{\tref}[1]{Theorem~\ref{#1}}
\newcommand{\dref}[1]{Diagram~\ref{#1}}

\def\Bbb#1{\mathbb#1}
\newcommand{\R}{\Bbb{R}}
\newcommand{\V}{\Bbb{V}}
\newcommand{\W}{\Bbb{W}}
\newcommand{\Q}{\Bbb{Q}}
\newcommand{\N}{\Bbb{N}}

\newcommand{\isim}{isometric immersion }
\newcommand{\isims}{isometric immersions }

\newcommand{\ra}{\rangle}
\newcommand{\la}{\langle}

\newcommand{\spa}{{\rm span \,}}
\newcommand{\im}{{\rm Im \,}}
\newcommand{\rk}{{\rm rank \,}}

\begin{document}

\title{Singular genuine rigidity}  
\author{Luis A. Florit and Felippe Guimar\~aes
\thanks{The authors were partially supported by CNPq-Brazil}}

\date{}
\maketitle

\begin{abstract}  
We extend the concept of genuine rigidity of submanifolds by allowing
mild singularities, mainly to obtain new global rigidity
results and unify the known ones.
\linebreak
As one of the consequences, we simultaneously extend and unify Sacksteder
and Dajczer-Gromoll theorems by showing that any compact $n$-dimensional
submanifold of ${\mathbb R}^{n+p}$ is singularly genuinely rigid in
${\mathbb R}^{n+q}$, for any $q < \min\{5,n\} - p$.
\linebreak
Unexpectedly, the singular theory becomes much
simpler and natural than the regular one, even though all technical
codimension assumptions, needed in the regular case, are removed.

\end{abstract}
\section{Introduction}

\hskip 0.6cm One of the fundamental problems in submanifold theory is
the (isometric) rigidity in space forms, i.e., whether an \isim of a
given Riemannian manifold is unique up to rigid motions. Satisfactory
solutions to the local version of the problem in low codimension were
obtained under certain nondegeneracy assumptions on the second
fundamental form, like the ones in \cite{Killing}, \cite{RigAll},
\cite{CDSnulidades}, \cite{DTComposition} and \cite{DFComposition}.
Recently, the concept of rigidity was extended to the one of genuine
rigidity in order to deal with deformations that arise as deformations of
submanifolds of larger dimension; see
\cite{DFGenuines}~and~\cite{FTConformal}. This reduction is important
since the difficulties in understanding rigidity aspects of submanifolds
grow together with the codimensions, not with the dimensions. This
concept also allowed to generalize and unify the papers mentioned above,
among others, by treating them under a common framework.

Global rigidity results are considerably more difficult to obtain. The
most important is the beautiful classical Sacksteder's theorem
\cite{SackstederBase}, which states that a compact Euclidean hypersurface
is rigid provided its set of totally geodesic points does not disconnect
the manifold. Outside the hypersurfaces realm there is only the
paper \cite{DGCompact}, where Dajczer and Gromoll showed that, along
each connected component of an open dense subset, any compact Euclidean
submanifold in codimension 2 is either genuinely rigid or a submanifold
of a special kind of deformable hypersurface. Although the authors did
not have the tools to justify it at the time, they had to allow certain
simple singularities in these hypersurfaces. The necessity to introduce
these singularities was justified recently in \cite{FGHonest}, and this
is precisely what motivated this work: to allow singularities in the
genuine rigidity theory, mainly with the double purpose of obtaining new
global results and unifying the known ones. In the process, we found out
that introducing these mild singularities is quite natural and
straightforward, even for local purposes, enabling us to substantially
simplify the theory. In fact, after completing this work, we regard the
presence of mild singularities in rigidity problems of submanifolds not
only as a necessary assumption to obtain global results but, more
importantly, as the natural setting for a deeper understanding
of the phenomena in an area where singularities rarely appear.

\medskip

In order to state our main results, let us introduce the key concepts.
We say that a pair of isometric immersions
$f:M^n\rightarrow\R^{n+p}$ and $\hat{f}: M^n \rightarrow\R^{n+q}$
{\it singularly extends isometrically}
when there are an embedding $j:M^n\hookrightarrow N^{n+s}$
into a manifold $N^{n+s}$ with $s>0$, and isometric maps
$F:N^{n+s}\to\R^{n+p}$ and $\hat{F}:N^{n+s}\to\R^{n+q}$
such that $f=F\circ j$ and $\hat f=\hat F\circ j$, with the set of points
where $F$ and $\hat{F}$ fail to be immersions (that may be empty)
contained in $j(M)$. In other words, the isometric extensions $F$ and
$\hat{F}$ in the following commutative diagram are allowed to be
singular, but only along $j(M)$:

\bigskip

\begin{picture}(110,84)\label{diag}
\put(152,31){$M^n$}
\put(206,31){$N^{n+s}$}
\put(410,31){$(1)$}
\put(239,62){$\R^{n+p}$}
\put(239,0){$\R^{n+q}$}
\put(192,59){$f$}
\put(192,4){${\hat f}$}
\put(237,46){${}_F$}
\put(236,26){${}_{\hat F}$}
\put(190,40.5){${}_j$}
\put(225,42){\vector(1,1){16}}
\put(225,28){\vector(1,-1){16}}
\put(174,44){\vector(3,1){60}}
\put(174,26){\vector(3,-1){60}}
\put(182,34){\vector(1,0){21}}
\put(182,36){\oval(7,4)[l]}
\end{picture}
\bigskip\medskip

\noindent An isometric immersion \immrw{n}{q} is a {\it strongly genuine
deformation} of a given isometric immersion \immr{n}{p}
if there is no open subset \mbox{$U \subset M^n$} along which the
restrictions $f |_U$ and $\hat{f} |_U$ singularly extend isometrically.
Accordingly, the isometric immersion $f$ is said to be {\it singularly
genuinely rigid} in $\R^{n+q}$ for a fixed integer $q$ if, for any given
isometric immersion \immrw{n}{q}, there is an open dense subset
$U\subset M^n$ such that $f|_{U}$ and $\hat{f}|_{U}$ singularly extend
isometrically.

More geometrically, an isometric deformation of a Euclidean submanifold
$M^n$ is strongly genuine if no open subset of $M^n$ is a submanifold of
a higher dimensional (possibly singular) isometrically deformable
submanifold, in such a way that the isometric deformation of the former
is induced by an isometric deformation of the latter, while (possibly)
including singularities along $M^n$. The key point here is that, since
all our extensions are ruled, the singularities that eventually appear
are quite mild and easy to understand, as it is classically done for the
classification of flat and ruled surfaces in $\R^3$.

\medskip

The following is our main global result. Recall that an immersion
$f$ is called {\it $D^d$-ruled}, or simply {\it $d$-ruled}, if
$D^d\subset TM$ is a rank $d$ totally geodesic distribution whose leaves
are mapped by $f$ to (open subsets of) $d$-dimensional affine subspaces.
Two immersions are said to be {\it mutually $d$-ruled} if they are
$D^d$-ruled with the same rulings $D^d$.

\begin{theorem}\po\label{CompactTheorem}
Let \immr{n}{p} and \immrw{n}{q} be isometric immersions of a compact
Riemannian manifold with $p+q<n$. Then, along each connected component of
an open dense subset of $M^n$, either $f$ and $\hat{f}$ singularly extend
isometrically, or $f$ and~$\hat{f}$ are mutually $d$-ruled, with
$d\geq n-p-q+3$.
\end{theorem}

In particular, for $p+q\leq 4$, \tref{CompactTheorem} easily unifies
Sacksteder and Dajczer-Gromoll Theorems in \cite{DGCompact} and
\cite{SackstederBase} cited above, states that the only way to
isometrically immerse a compact Euclidean hypersurface in codimension 3
is through compositions (which in turn were classified in \cite{cfs} and
\cite{FGHonest}), and provides a global version of the main
result~in~\cite{DFSbrana}:

\begin{corollary}\po \label{DGGeneral}
Any compact isometrically immersed submanifold $M^n$ of $\R^{n+p}$ is
singularly genuinely rigid in $\R^{n+q}$ for $q < min\{5,n\}-p$.
\end{corollary}

 From \tref{CompactTheorem} we get the following topological criteria for
singular genuine rigidity in line with the rigidity question proposed by
M. Gromov in \cite{G} p.259 and answered in \cite{DGCompact} (and thus
also in \cref{DGGeneral}), without any {\it a priori} assumption on the
codimensions:

\begin{corollary}\po \label{theoExtra}
Let $M^n$ be a compact manifold whose $k$-th Pontrjagin class
satisfies that $[p_k] \neq 0$ for some $k> \frac{3}{4}(p+q-3)$. Then, any
analytic immersion \immr{n}{p} (with the induced metric) is singularly
genuinely rigid in $\R^{n+q}$ in the $C^{\infty}$-category.
\end{corollary}

Our global results are based on a local analysis whose main tool is the
bilinear form that we construct next. Consider a pair of isometric
immersions \immr{n}{p} and \immrw{n}{q}. Let
$$
\tau\colon L^\ell\subset T^{\perp}_fM\rightarrow
\hat{L}^\ell\subset T^{\perp}_{\hat{f}}M
$$
be a vector bundle isometry and suppose that it preserves the second
fundamental forms and the normal connections restricted to the rank
$\ell$ vector normal subbundles $L^\ell$~and~$\hat L^\ell$.
Equivalently, its natural extension
$\bar\tau=Id\,\oplus\tau\colon TM\oplus L^\ell\to TM\oplus\hat L^\ell$
is a parallel bundle isometry. Let
$\phi_\tau\colon TM \times(TM\oplus L^\ell)\rightarrow
L^\perp\times\hat L^\perp$
be the flat bilinear form given by
$$
\phi_\tau(X,v) = \left((\tilde{\nabla}_X v)_{L^{\perp}},
(\tilde{\nabla}_X \bar \tau v)_{\hat{L}^{\perp}}\right),
\ \ \ X \in TM, \ v \in TM\oplus L^\ell,
$$
where $\tilde{\nabla}$ stands for the connection in Euclidean space and
$L^\perp\times\hat L^\perp$ is endowed with the semi-Riemannian metric
$\langle\,,\rangle=\langle\,,\rangle|_{L^\perp}
-\langle\,,\rangle|_{\hat L^\perp}$.
A subset $S\subset L^\perp\oplus\hat L^\perp$ is called {\it null}
if $\langle\eta,\xi\rangle = 0$ for all $\eta,\xi \in S$.

\medskip

In order to present our local statements we need to extend the concept of
$D^d$-ruled for arbitrary distributions $D^d\subset TM$. In this case, we
say that $f$ is {\it $\bar D^d$-ruled} if, for each $p\in M^n$, there is
a totally geodesic submanifold of $M^n$ tangent to $D(p)$ at
$p$ which is mapped by $f$ to an (open subset of) a $d$-dimensional
affine subspace. Observe that such an $f$ is of course also $d$-ruled as
before, but usually with bigger rulings when $D^d$ is not totally
geodesic.

We can now state our main local result, which applies even to $\ell=0$ and
$\tau = 0$.

\begin{theorem}\po \label{MainTheorem}
Let \immrw{n}{q} be a strongly genuine deformation of \immr{n}{p} and
$\tau:L^\ell \subset T^{\perp}_fM \rightarrow \hat L^\ell \subset
T^{\perp}_{\hat{f}}M$ a parallel vector bundle isometry that preserves
second fundamental forms. Let $D \subset TM \oplus L^\ell$ be a subbundle
such that $\phi_\tau(TM,D)$ is a null subset. Then $D \subset TM$ and,
along each connected component of an open dense subset of $M^n$, $f$ and
$\hat{f}$ are mutually $\bar D$-ruled.
\end{theorem}

The usefulness of \tref{MainTheorem} relies on the fact that it deals
with easily to construct null subsets instead of nullity distributions of
flat bilinear forms. A good example of an application of this fact is the
following singular version of Theorem~1 in \cite{DFGenuines} removing the
technical assumption on the codimensions.
Recall that $Y \in T_xM$ is a {\it regular element} of
$\phi_\tau$ at $x$ if $\rk(\phi_\tau^Y)=i(\phi_\tau)(x)$, where
$\phi_\tau^Y = \phi_\tau(Y,\,\cdot\,)$ and
$$
i(\phi_\tau)(x):= \max \{ \text{rank}(\phi_\tau^X)\colon X \in T_xM \}.
$$
Denote by $RE(\phi_\tau) \subset TM$ the open dense subset of
regular elements of $\phi_\tau$. Using a well-known property of flat
bilinear forms we immediately conclude from \tref{MainTheorem}:

\begin{corollary}\po \label{CorMain} Under the assumptions of
\tref{MainTheorem}, along each connected component of an open dense
subset of $M^n$, $i(\phi_\tau)$ is constant and $f$~and~$\hat{f}$ are
mutually $\bar D_Y^d$-ruled for any smooth vector field
$Y \in RE (\phi_\tau)$, where $D_Y^d := \ker(\phi^{Y}_{\tau})\subset TM$.
In particular, $f$~and~$\hat{f}$ are
mutually $d$-ruled with $d=n+\ell-i(\phi_\tau)\geq n-p-q+3\ell$.
\end{corollary}

As it is clear from the statements, the rulings in the above are larger
and easier to compute than the ones in the main result in
\cite{DFGenuines}. The bundles obtained in this work are also better
suited for certain global applications.


\medskip

By allowing singular extensions we recover all the corollaries in
\cite{DFGenuines}, even without the technical restrictions on the
codimensions required there. For example, from \cref{CorMain} we conclude
the following extension of Corollary~5 in \cite{DFGenuines}.

\begin{corollary}\po
Any isometrically immersed submanifold $M^n$ of $\R^{n+p}$ with
positive Ricci curvature is singularly genuinely rigid in $\R^{n+q}$,
for every $q<n-p$.
\end{corollary}


As we will see, the proof of the local Theorem $\ref{MainTheorem}$ works
for any simply connected space form. Moreover, the global
\tref{CompactTheorem} and \cref{DGGeneral} still hold for
complete submanifolds, even if the ambient space is the hyperbolic space,
as long as one of the immersions is bounded. For complete submanifolds in
the round sphere we show:

\begin{theorem}\po \label{GloRes}
Let $f: M^n \rightarrow \mathbb{S}^{n+p}$ and $\hat{f}: M^n \rightarrow
\mathbb{S}^{n+q}$ be isometric immersions of a complete submanifold with
$p+q<n-\mu_n$. Then, along each connected
component of an open dense subset of $M^n$, either $f$ and $\hat{f}$
singularly extend isometrically, or $f$ and $\hat{f}$ are mutually
$d$-ruled, with $d\geq n-p-q+3$.
\end{theorem}

In the above statement $\mu_n$ is defined as
$\mu_n = \max\{ k\, : \, \rho(n-k) \geq k + 1\}$,
where $\rho(m)-1$ is the maximum number
of pointwise linearly-independent vector fields on $\mathbb{S}^{m-1}$
and is given by $\rho((odd)2^{4d+b}) = 8d + 2^b$,
for any nonnegative integer $d$ and $b \in \{0,1,2,3\}$. Some values of
$\mu_n$ are: $\mu_n = n - {\rm (highest\ power\ of\ }2 \leq n$) for
$n\leq 24$, $\mu_n \leq 8d-1$ for $n<16^d$ and $\mu_{2^d} = 0$.

\medskip

 From \tref{GloRes} we obtain the corresponding version of
\cref{DGGeneral} for complete submanifolds in the sphere:

\begin{corollary}\po \label{corEsf}
Any complete isometrically immersed submanifold $M^n$ of
$\mathbb{S}^{n+p}$ is singularly genuinely rigid in
$\mathbb{S}^{n+q}$ for $q \leq 3 - p$ if $4 \leq n \leq 7$, or $q\leq
4-p$ if $n \geq 8$.
\end{corollary}

The paper is organized as follows. In \sref{mbm} we first provide the
basic properties of the bilinear form $\phi_\tau$, and then we show how
it can be used to obtain regular and singular isometric extensions,
which is all that is needed to prove our local results. \sref{compos} is
devoted to revisit the theory of compositions using $\phi_\tau$. As an
application we show that, generically, $(n-1)$-ruled submanifolds are
compositions. In \sref{ga} we prove \tref{CompactTheorem}, and
\sref{sphereCase} is dedicated to the proof of \tref{GloRes}.

\setcounter{equation}{1}

\section{The flat bilinear form \texorpdfstring{$\phi_\tau$}{phi}}\label{mbm}

\hskip 0.6cm In this section we study some properties of the
bilinear form $\phi_\tau$, which was introduced in \cite{DFGenuines} but
not used in its full strength. We will see that it is a powerful tool
to deal with isometric rigidity problems.

\medskip

Consider two isometric immersions \immr{n}{p} and \immrw{n}{q} with second
fundamental forms $\alpha$ and $\hat{\alpha}$ and normal connections
$\nabla^\perp$ and $\hat\nabla^\perp$ defined on their normal bundles
$T^{\perp}_fM$ and $T^{\perp}_{\hat f}M$, respectively. Endow
$T^{\perp}_fM \times T^{\perp}_{\hat{f}}M$ with its natural
semi-Riemannian metric $\langle\, ,\rangle$ of type $(p,q)$
and compatible connection $\bar \nabla$,
$$
\langle (\xi,\hat{\xi}),(\eta,\hat{\eta})\rangle
=\langle\xi,\eta\rangle_{T^{\perp}_fM}
-\langle\hat{\xi},\hat{\eta}\rangle_{T^{\perp}_{\hat{f}}M},\ \ \ \ \
\bar \nabla_X(\xi,\hat \xi)
=\left(\nabla^\perp_X\xi,\hat\nabla^\perp_X\hat\xi\right),
$$
for $\xi,\eta\in T^{\perp}_fM$, $\hat\xi,\hat\eta\in T^{\perp}_{\hat f}M$,
and $X\in TM$. By the Gauss equation, the symmetric bilinear form
$\beta=\alpha \oplus \hat{\alpha}:
TM\times TM\rightarrow T^{\perp}_fM\times T^{\perp}_{\hat{f}}M$
is {\it flat}, that is,
$$
\langle \beta(X,Y), \beta(Z,T)\rangle = \langle \beta(X,T),
\beta(Z,Y)\rangle, \, \forall X,Y, Z, T \in TM.
$$
The concept of flat bilinear forms was introduced by Moore in
\cite{Moore} to study isometric immersions of the round sphere in
Euclidean space in low codimension, and was used afterwards in several
papers about isometric rigidity, even implicitly, following
a remark also in \cite{Moore}. For example, it can be
used to prove the classical Beez-Killing theorem in~\cite{Killing}, in
which case the objective is to show that $\im(\beta)$ is everywhere a
null set. Notice that flatness makes sense even for nonsymmetric bilinear
forms.

Outside the realm of hypersurfaces it is important to obtain information
about the normal connections too, so a different (nonsymmetric) flat
bilinear form is needed. Yet, unexpectedly and in contrast to the
strongest known local rigidity results, we will not make use
of {\it a priori} nullity estimates like the one in Theorem 3 in
\cite{DFComposition} in ours since we will not deal with nullity spaces.
In particular, this will allow us to get rid of the usual technical
constraints on the codimensions.

\medskip

Throughout this work,
$$
\tau: L^\ell \subset T^{\perp}_fM \rightarrow
\hat{L}^\ell \subset T^{\perp}_{\hat{f}}M
$$
will denote a vector bundle isometry that preserves the induced second
fundamental forms and normal connections in the rank $\ell$ normal
subbundles $L$ and $\hat{L}$. That is,
$\tau\circ\alpha_L=\hat\alpha_{\hat L}$, and
$\tau (\nabla_X^{\perp} \xi)_{L}=(\hat{\nabla}_X^\perp\tau\xi)_{\hat{L}}$
for every $X \in TM$, $\xi \in L$,
where we represent the orthogonal projections onto $L$ and $\hat{L}$ with
the corresponding subindexes. Equivalently, its natural extension
$$
\overline{\tau}=Id\oplus\tau:TM\oplus L\rightarrow TM\oplus\hat{L},
$$
is a parallel vector bundle isometry. Let
$\phi_\tau: TM \times (TM\oplus L)\rightarrow
L^{\perp}\times\hat{L}^{\perp}\subset
T_f^{\perp}M\times T_{\hat{f}}^{\perp}M$
be the bilinear form defined as
$$
\phi_\tau(X,v)=\left((\tilde{\nabla}_X v)_{L^{\perp}},
(\tilde{\nabla}_X \overline{\tau}v)_{\hat{L}^{\perp}}\right),
$$
where $\tilde{\nabla}$ denotes the connection of the Euclidean ambient
spaces.
Notice that, if $\ell=0$, then $\tau=0$ and $\phi_0 = \beta$.
On $L^{\perp} \times \hat{L}^{\perp}$ we will always consider the
semi-Riemannian metric and compatible connection induced from the ones
in $T^{\perp}_fM \times T^{\perp}_{\hat{f}}M$, still denoted by
$\langle\, ,\rangle$ and $\bar\nabla$, respectively.

\medskip

The main two properties of $\phi_\tau$ are given by the following.
\begin{proposition}\po\label{mainprops}
$\phi_\tau$ is a flat Codazzi tensor.
\end{proposition}
\proof
For $X,Y \in TM$ and $v,w \in TM \oplus L$, using that $\bar \tau$ is
parallel we get
\begin{equation*}
\begin{split}
\langle\phi_\tau(X,v),\phi_\tau(Y,w) \rangle
&=\langle(\tilde{\nabla}_Xv)_{L^{\perp}},(\tilde{\nabla}_Yw)_{L^{\perp}}\rangle
-\langle (\tilde{\nabla}_X\overline{\tau}v)_{\hat{L}^{\perp}},
(\tilde{\nabla}_Y\overline{\tau}w)_{\hat{L}^{\perp}} \rangle \\
&=\langle\tilde{\nabla}_Xv,\tilde{\nabla}_Yw\rangle
-\langle\tilde{\nabla}_X\overline{\tau}v,\tilde{\nabla}_Y\overline{\tau}w \rangle \\
&=-\langle v,\tilde{\nabla}_Y\tilde{\nabla}_Xw\rangle
+\langle\overline{\tau}v,\tilde{\nabla}_Y\tilde{\nabla}_X\overline{\tau}w\rangle\\
&=-\langle v,\tilde{\nabla}_X\tilde{\nabla}_Yw\rangle
+\langle \overline{\tau}v,\tilde{\nabla}_X \tilde{\nabla}_Y\overline{\tau}w \rangle \\
&=\langle \phi_\tau(Y,v),\phi_\tau(X,w) \rangle.
\end{split}
\end{equation*}
The very same approach shows that $\phi_\tau$ is a Codazzi tensor, that
is,
$$
(\bar \nabla_X \phi_\tau)(Y,v):=
\bar \nabla_X \phi_\tau(Y,v)-\phi_\tau(\nabla_XY,v)-\phi_\tau(Y,(\tilde\nabla_X v)_{TM\oplus L}) = (\bar \nabla_Y \phi_\tau)(X,v),
$$
so we left the computation to the reader.
\qed
\medskip

Denote the left nullity space of $\phi_\tau$ by $\Delta_\tau$ and its
dimension by $\nu_\tau$, i.e.,
\begin{equation}\label{deltatau}
\Delta_\tau:=\{X\in TM: \phi_\tau(X,\,\cdot\,)=0\}, \ \ \ \
\nu_\tau:=\dim\Delta_\tau,
\end{equation}
and let $U\subset M^n$ be a
connected component of an open dense subset where $\nu_\tau$ is
locally constant. Since $\phi_\tau$ is a Codazzi tensor,
$\Delta_\tau$ is a smooth integrable distribution on $U$.

\begin{corollary}\po\label{par}
The space $\im(\phi_\tau)^\perp$ is
parallel along the leaves of $\Delta_\tau$ in $U$.
In particular, both the nullity space and the light cone bundle of
$\la\,,\,\ra|_{\im(\phi_\tau)^\perp}$ are smooth and parallel along
these leaves on any open subset $U'\subset U$ where they have constant
dimension.
\end{corollary}
\proof
The parallelism of $\im(\phi_\tau)^\perp$
is a consequence of the fact that $\phi_\tau$ is a Codazzi tensor, since
$\bar \nabla_X (\phi_\tau(Y,v))=\phi_\tau([X,Y],v)
+\phi_\tau(Y,(\tilde\nabla_X v)_{TM\oplus L})\in\spa\im(\phi_\tau)$
for every $X\in\Delta_\tau$, $Y\in TM$, $v\in TM\oplus L$.
The last assertion follows from the compatibility of $\bar\nabla$
with respect to $\la\,,\,\ra$.
\qed
\medskip

Observe that such a $\tau$ as above arises naturally when $f$ and
$\hat f$ singularly extend isometrically. Indeed, with the
notations in \dref{diag} in the Introduction, if $F$ and $\hat F$ are
regular we just take $L^s=F_*(T_j^\perp M)$,
$\hat L^s=\hat F_*(T_j^\perp M)$, and $\tau(F_*(\xi))=\hat F_*(\xi)$.
If they are not regular, at least locally we can consider a sequence of
submanifolds $j_k\colon M^n_k\to N^{n+s}\setminus j(M^n)$ smoothly
converging to $j$ as $k\to\infty$, and then take $L^s$ and $\hat L^s$ as
accumulation of $F_*(T_{j_k}^\perp M_k)$ and
$\hat F_*(T_{j_k}^\perp M_k)$, respectively. In particular, we have:

\begin{lemma}\po\label{max}
The metric induced on
$\im(\beta)^\perp\subset T_f^\perp M\oplus T_{\hat f}^\perp M$
is almost everywhere definite if and only if $\tau=0$ is locally the only
vector bundle isometry preserving second fundamental forms. In this
situation, $\hat f$ is a strongly genuine deformation of $f$.
\end{lemma}
\proof
Both conditions are clearly equivalent to the non-existence of unit
vector fields $\xi\in T^\perp_f M$ and $\hat \xi\in T^\perp_{\hat f} M$
defined on some open subset $U\subset M^n$ such that $A_\xi=\hat A_{\hat \xi}$.
\qed

\subsection{The form \texorpdfstring{$\phi_\tau$}{phi} and genuine rigidity} 

\hskip 0.6cm In general, we show that a pair of \isims $\{f,\hat f\}$ as
above is genuine, i.e., each one is a (regular) genuine deformation of
the other, by explicitly constructing, locally almost everywhere,
isometric immersions $j:M^n\to N^{n+s}$, $F:N^{n+s}\to\R^{n+p}$
and $\hat F:N^{n+s} \to \R^{n+q}$ as in \dref{diag}, that is,
satisfying $f=F\circ j$ and $\hat f=\hat F\circ j$.
Usually, we also require $F$ and $\hat F$ to be ruled extensions of
$f$ and $\hat f$ since a genuine pair must be mutually ruled by the main
result in \cite{DFGenuines}. Since in this paper we work with singular
extensions, the ruled ones have the additional advantage that their
singularities are quite easy to characterize and deal with.

\medskip

In order to build ruled extensions of $f$ and $\hat f$, choose any
smooth rank $s$ subbundle
$\Lambda \subset TM \oplus L$, and define the maps $F = F_{\Lambda,f}:
\Lambda \rightarrow \R^{n+p}$ and
$\hat{F} = F_{\Lambda, \hat{f}}:\Lambda \rightarrow \R^{n+q}$ as
\begin{equation} \label{ExtNat}
F(v) = f(p) + v, \ \ \ \hat{F}(v) = \hat{f}(p) +
\overline{\tau}v, \ \ \ v \in \Lambda_p,\ \ p\in M^n.
\end{equation}
One of the main reasons that make the form $\phi_\tau$ useful in any
flavour of genuine rigidity is that it gives the precise condition that
guarantees that these two maps are isometric:
\begin{proposition}\po\label{iso}
The maps $F$ and $\hat F$ in \eqref{ExtNat} are isometric if and only
if $\phi_\tau(TM,\Lambda)$ is a null set.
\end{proposition}
\proof
It follows easily from the fact that $\overline{\tau}$ is
parallel since, for every smooth local section $v$ of $\Lambda$ and $Z
\in TM$, we have that
$\|(F\circ v)_*Z\|^2-\|(\hat{F}\circ v)_*Z\|^2
=\|Z+\tilde{\nabla}_Zv\|^2-\|Z+\tilde{\nabla}_Z\overline{\tau}v\|^2
=\|(\tilde{\nabla}_Z v)_{L^{\perp}} \|^2
-\|(\tilde{\nabla}_Z \overline{\tau}v)_{\hat{L}^{\perp}} \|^2
=\la\phi_\tau(Z,v),\phi_\tau(Z,v)\ra$.
\qed

\medskip

In particular, if in addition $\Lambda \cap TM = 0$, both maps are
immersions in a neighborhood $N^{n+s}$ of the $0$-section of $\Lambda$,
and thus induce the same Riemannian metric on $N^{n+s}$. Therefore $F$
and $\hat{F}$ are (regular) isometric ruled extensions of $f$ and
$\hat{f}$. Similarly, if $\Lambda\not\subset TM$, along each open subset
$U\subset M^n$ where the subspaces $\Lambda'=\Lambda\cap(\Lambda\cap
TM)^\perp$ have locally constant dimension $s'>0$, we have that the
restrictions $F'|_{\Lambda'}=F_{\Lambda',f|_U}$ and
$\hat F'|_{\Lambda'}=F_{\Lambda',\hat f|_U}$
also give (regular) isometric ruled extensions of $f|_U$ and $\hat f|_U$
defined in a neighborhood $N^{n+s'}$ of the $0$-section of $\Lambda'$
along $U$.

\medskip

We proceed to characterize singular ruled extensions, that occur above
when $\Lambda \subset TM$. We say that $F = F_{\Lambda,f}$ in
\eqref{ExtNat} is a {\it singular extension of} $f$ if it is an immersion
in some open neighborhood of the $0$-section of $\Lambda$, except of
course at the $0$-section itself. We say that $F$ {\it nowhere induces
a singular extension of} $f$ if, for every open subset $U \subset M^n$
and every subbundle $\Lambda' \neq 0$ of $\Lambda|_U$, the restriction of
$F$ to $\Lambda'$ is not a singular extension of $f|_U$. We show next
that $F$ nowhere induces a singular extension of $f$ only when the latter
is $\bar\Lambda$-ruled.

\begin{proposition}\po\label{singD}
Let \immr{n}{p} be an isometric immersion and $\Lambda \subset TM$
a smooth distribution. Then, $F_{\Lambda,f}$ nowhere induces a singular
extension of $f$ if and only if $f$ is $\bar\Lambda$-ruled along
each connected component of an open dense subset of $M^n$.
\end{proposition}
\proof
Clearly, it is enough to give a proof for the direct statement and for a
rank one distribution, i.e., $\Lambda = \spa\{X\}$ for some nonvanishing
vector field $X$ on $M^n$.
Consider the map $F: \Lambda \cong M^n \times \R \rightarrow \R^{n+p}$
given by (\ref{ExtNat}), that is, $F(p,t)=f(p)+tX(p)$. This map will be a
singular extension in some open neighborhood of $p \in M^n$ if and only
if it is an immersion in a neighborhood of $(p,0)$, except at the points
in $M^n \times \{0\}$. Therefore, for all $p \in M^n$ there exists a
sequence $(p_m, t_m) \rightarrow (p,0)$, with $t_m \neq 0$, such that
rank$(F_{*(p_m, t_m)}) = n$. Define the tensors $K(Z) = \nabla_Z X$ and
$H_t(Z) = Z + tK(Z)$ for $Z \in TM$. Thus, there is $Y_m \in T_{p_m}M$
such that $ F_{*(p_m, t_m)} Y_m = X(p_m)$, since $H_t \rightarrow Id$ as
$t \rightarrow 0$ and
$$
F_*\partial_t=X,\ \ \ \ F_* Z = H_t(Z)+t\alpha(X,Z),
\ \ \forall\, Z \in TM.
$$

Let $S_X$ be the $K$-invariant subspace generated by $X$,
$$
S_X = \spa\{X, K(X), K^2 (X), K^3(X), ... \}.
$$
Observe that the equality
$F_{*(p_m,t_m)} Y_m = X(p_m)$ is equivalent to $H_{t_m}Y_m = X(p_m)$ and
$\alpha(X(p_m),Y_m)=0$. In particular, if $t_m$ is sufficiently small,
\begin{equation}\label{propSingEq}
\alpha(X(p_m),H_{t_m}^{-1}(X(p_m))) = 0
\end{equation}
and $\lim_{m \rightarrow \infty} H^{-1}_{t_m}(X(p_m)) = X(p)$. Consider a
precompact open neighborhood $U \subset M^n$ of $p$, so $\| \alpha
\| < c$ and $\| K \| < c$ for some constant $c>1$. Hence for $t
\in I = (-\frac{1}{c^2},\frac{1}{c^2}) $ we have that $H_t$ is invertible
on $U$, and
$$
H_t^{-1} = \sum_{i\geq0} (-t)^i K^i,
$$
since $H_t(\sum_{i=0}^N (-t)^i K^i) = Id - (-t)^{N+1}K^{N+1}$.

We claim that $\alpha(X, S_X) = 0$ along $M^n$. Assume otherwise,
define $j:=\min\{k\in\N:\alpha(X(q),K^k(X(q)))\neq0,\ q\in M^n\}$
and take $p\in M^n$ such that $\alpha(X(p),K^j (X(p)))\neq0$.
By~\eqref{propSingEq} we obtain that
$$\sum_{i\geq j} (-t_m)^i \alpha(X(p_m), K^i (X(p_m))) =0.
$$
Dividing the above by $t^j_m$ and taking $m\rightarrow\infty$ we conclude
that $\alpha(X(p),K^j (X(p)))=0$, which is a contradiction.

Now, since $\alpha(X ,S_X) = 0$ on $M^n$, for any $t \in I$ and $p \in U$
we get $F_{*(p,t)}(H^{-1}_t(X)) = X$ since $H^{-1}_t(X) \in S_X$.
It follows that rank$(F_*) = n$ in all $U \times I$, and therefore
$F(U\times I) = f(U)$. Hence a segment of the line generated by $X$ is
contained in $f(U)$.
\qed
\medskip

We are now able to prove our main local result.

\medskip

\noindent {\it Proof of \tref{MainTheorem}.} Locally, if
$D\not\subset TM$ along some open set $U$ then we have regular isometric
extensions of $f|_U$ and $\hat{f}|_U$ by extending them as in
\eqref{ExtNat} along any subbundle $\Lambda \subset D$ such that
$D = (D \cap TM) \oplus \Lambda$. Hence, $D\subset TM$ and by
\pref{singD} we conclude that $f$ and $\hat{f}$ are mutually
$\bar D$-ruled almost everywhere.
\qed
\medskip

The following lemma due to Moore \cite{Moore}
immediately gives \cref{CorMain} by applying \tref{MainTheorem}
to $\tau$ and $D^d = \ker(\phi_\tau^Y)$, since it tells us that
$\phi_\tau(TM,\ker(\phi_\tau^{Y}))$ is null.

\begin{lemma}\po\label{mr}
Let $\varphi:\V\times\V'\to \W$ be a flat bilinear form, and set
$\varphi^X=\varphi(X,\,\cdot\,)$. Then,
$$
\varphi(\V,\ker(\varphi^X))\subset
\im(\varphi^X)\cap \im(\varphi^X)^\perp,\ \ \
\forall\,X\in RE(\varphi).
$$
In particular, if the inner product in $\W$ is definite, we have
that $\, \ker(\varphi^X)={\cal N}(\varphi)$ for all $X\in RE(\varphi)$,
where ${\cal N}(\varphi):=\{w\in\V':\varphi(\,\cdot\,,w)=0\}$
is the (right) nullity of $\varphi$.
\end{lemma}

%

\begin{remark}\po\label{easy}{\rm
While \cref{CorMain} with its estimate $d\geq n-p-q+3\ell$ is immediate
from \lref{mr}, the corresponding regular result, Theorem 14 in
\cite{DFGenuines}, requires several pages just to give a proof of the
estimate on $d$. In addition, it uses the very long and technical Theorem
3 in \cite{DFComposition}, and therefore it is only valid for
$\min\{p,q\}\leq 5$; see \cite{count}. The simplifications
gained with the singular theory reside in the fact that, while here we
use \lref{mr} to easily obtain null subsets, the main results in
\cite{DFGenuines} require the computation of estimates of ranks of
several bundles and nullities of trickily constructed bilinear forms.}
\end{remark}

\begin{remark}\po\label{structureRmk}{\rm
In several applications we have that $D={\cal N}(\alpha_{L^\perp})\cap
{\cal N}(\hat{\alpha}_{\hat{L}^{\perp}})$ even in the singular case.
For example, this is the case if $d=n-p-q+3\ell$ in \cref{CorMain}, or if
$\ell=\min\{p,q\}$, or if one of the codimensions is low enough. In this
situation, $L_D\subset L$, $\tau|_{L_D}$ is also
parallel and preserves second fundamental forms, and therefore we recover
the structure of the normal bundles in Theorem 1 in \cite{DFGenuines};
see \lref{l} below.}
\end{remark}

\section{Compositions revisited through \texorpdfstring{$\phi_\tau$}{phi}}\label{compos} 
\hskip 0.6cm In this section we revisit the theory of compositions using
the form $\phi_\tau$.

\medskip

Let $f:M^n\to\R^{n+p}$ be an isometric immersion of a simply connected
Riemannian manifold $M^n$ with second fundamental form $\alpha$, and
$L\subset T^\perp_fM$ a rank $\ell$ normal subbundle.
Define the bilinear form
$$
\phi_{L^\perp}:TM\times (TM\oplus L)\to L^\perp, \ \ \ \
\phi_{L^\perp}(Z,v)=(\tilde\nabla_Zv)_{L^\perp}.
$$
We can build another \isim of $M^n$ using $\phi_{L^\perp}$ when it is flat:
\begin{proposition}\po\label{exists}
The bilinear form $\phi_{L^\perp}$ is flat if and only if there exists an
isometric immersion $\hat f: M^n\to\R^{n+\ell}$ and a parallel vector
bundle isometry $\sigma:L\to T^\perp_{\hat f}M$ such that the second
fundamental form of $\hat f$ is $\hat\alpha=\sigma\circ\alpha_L$.
In this case, $\phi_{L^\perp}=\phi_\sigma$.
\end{proposition}

\proof
By projecting the fundamental equations of $f$ onto $L$ we easily see
that flatness of $\phi_{L^\perp}$ is equivalent to the fact that the pair
$(\alpha_L,(\nabla^\perp)|_L)$ satisfies the fundamental equations
of Euclidean submanifolds. Indeed, flatness of $\phi_{L^\perp}$ with the
four vectors in $TM$ is equivalent to Gauss equation, with three vectors
in $TM$ and one in $L$ we get Codazzi equation, while two vectors in $TM$
and two in $L$ recovers Ricci equation.
\qed
\medskip

The following is a reinterpretation of Proposition 8 in
\cite{DFComposition}, which is the main tool to construct compositions.
Recall that, for $f$ and $\hat f$ as in the previous section, we say
that~$f$ is a {\it (regular) composition} of $\hat f$ when they extend
isometrically as in \dref{diag} with $s=q$. In this case, $\hat F$ is a
local isometry and thus, if $\hat f$ is an embedding, there is an open
neighborhood $U\subset N^{n+q}$ of $j(M^n)$ and an isometric immersion
$h=F|_U\circ(\hat F|_U)^{-1}:W\subset\R^{n+q}\to\R^{n+p}$ of the open
subset $W=\hat F(U)\supset f(M)$ satisfying $f=h\circ\hat f$. Recall also
that $i(\phi_{L^\perp})(x):=\max\{\rk(\phi_{L^\perp}^X):X\in T_xM\}$.

\begin{proposition}\po\label{comp2}
Suppose that $\phi_{L^\perp}$ is flat, and let $\hat f$ be given by
\pref{exists}. \linebreak If $i(\phi_{L^\perp})$ is constant and
$i(\phi_{L^\perp})=i(\alpha_{L^\perp})$, then $f$ is a composition of
$\hat f$.
\end{proposition}
\proof
Observe that, since $\hat L=T^\perp_{\hat f}M$ and
$\phi_{L^\perp}=\phi_\sigma$, the image of both $\phi_{L^\perp}$ and
$\alpha_{L^\perp}=\phi_{L^\perp}|_{TM\times TM}$ are Riemannian. Thus, by
\lref{mr} we have that
$\ker(\phi_{L^\perp}^X)={\cal N}(\phi_{L^\perp})$ and
$\ker(\alpha_{L^\perp}^X)={\cal N}(\alpha_{L^\perp})$,
for every $X\in RE(\phi_{L^\perp}) \cap RE(\alpha_{L^\perp})$.
The result follows from \pref{iso} taking the rank $\ell$
subbundle $\Lambda={\cal N}(\phi_{L^\perp})\cap{\cal
N}(\alpha_{L^\perp})^\perp$, which is transversal to $M^n$.
\qed

\begin{remark}\po
{\rm By allowing singular flat extensions $h$, \pref{singD} provides a
singular version of \pref{comp2}: if there is a rank $\ell$ subbundle
$\Lambda\subset{\cal N}(\phi_{L^\perp})$ such that $F_{\Lambda,f}$ is an
immersion near the 0-section of $\Lambda$, except possibly along it, then
$f$ is locally almost everywhere a {\it singular composition}
$f=h\circ \hat f$, where $h$ is an immersion except (possibly) along
$\hat f(M)$.}
\end{remark}

\begin{remark}\po
{\rm By \tref{MainTheorem} applied to $\tau=\sigma$, if $\hat f$ in
\pref{exists} is a strongly genuine deformation of $f$, then they must be
at least mutually $(n-p+2\ell)-$ruled, and by \pref{comp2},
$i(\alpha_{L^\perp})<i(\phi_{L^\perp})\leq p-\ell$.}
\end{remark}


\subsection{The \texorpdfstring{$(n-1)$}{(n-1)}-ruled case}\label{n1} 
\hskip 0.6cm As an application of the above, here we study
general $(n-1)$-ruled $n$-dimensional Euclidean submanifolds. We show
that such a submanifold is locally a composition if its codimension is
bigger than the rank of its curvature operator. Although this fact has
independent interest, it will be used to prove \cref{DGGeneral}.
\medskip


Until the end of this section $X,Y$ will denote vectors in a totally
geodesic distribution $D\subset TM$, and $Z\in TM$.
\begin{lemma}\po\label{l}
If $f$ is $D$-ruled, then the normal subbundle
\begin{equation}\label{ld}
L_D :=\spa\alpha(TM,D)\subset T_f^\perp M
\end{equation}
is parallel along $D$ on any open subset $V$ where $\ell_D:=\dim L_D$
is constant.
\end{lemma}
\proof
Since $D$ is totally geodesic, the lemma follows from Codazzi equation since
$$
\nabla^\perp_X\alpha(Z,Y)=-\alpha(\nabla_XZ,Y)-\alpha(Z,\nabla_XY)
-\alpha(\nabla_ZX,Y)-\alpha(X,\nabla_ZY)\in L_D.
\qed
$$

In particular, if rank $D=n-1$, $L_D\subset L$ and $L$ is also parallel
along $D$, our form $\phi_{L^\perp}$ is flat since
$\phi_{L^\perp}(D^{n-1},TM\oplus L)=0$. Therefore \pref{exists} gives:

\begin{corollary}\po\label{ruled}
Suppose $f$ is $D^{n-1}$-ruled and $L\subset T_f^\perp M$ is a rank
$\ell$ normal subbundle parallel along $D^{n-1}$ such that
$L_D\subset L$. Then, there is a $D^{n-1}$-ruled isometric immersion
$$
\hat f:M^n\to\R^{n+\ell}
$$
and a parallel vector bundle isometry
$\sigma:L\to T^\perp_{\hat f}M$ such that the second fundamental form of
$\hat f$ is $\hat \alpha=\sigma\circ\alpha_L$. In particular, taking $V$
in \lref{l} simply connected, there exists a $D^{n-1}$-ruled isometric
immersion $f_D:V\subset M^n\to\R^{n+\ell_D}$ with
$\hat\alpha=\sigma\circ\alpha_{L_D}$.
\end{corollary}

Notice that, if rank $D =n-1$, $\ell_D$ is intrinsic since it agrees
with the rank of the curvature operator of $M^n$.
Our purpose is to show that $f$ is locally a composition of $f_D$.

\begin{proposition}\po\label{ruled2}
Under the assumptions of \cref{ruled}, if $\im \alpha(x)\not\subset L(x)$
for some $x\in M^n$, then $f$ is a composition of $\hat f$ near $x$, that
is, there is a neighborhood $U$ of $x$ and an isometric immersion
$h:W\subset\R^{n+\ell}\to\R^{n+p}$ of an open set $W$ containing $f(U)$
such that $f=h\circ \hat f$ on $U$.
\end{proposition}
\proof
The hypothesis is equivalent to the existence of an orthogonal
decomposition
$$
T_f^\perp V = L \oplus \Sigma \oplus N,
$$
on an open neighborhood $V$ of $x$, where $\Sigma$ is a line bundle and
$N=\{\eta\in L^\perp:A_\eta=0\}$. We proceed by induction on the
codimension $p\geq\ell+1$ of $f$. For $p=\ell+1$ we get that $f$
is a composition of $\hat f$ near $x$ by \pref{comp2}
since $1\leq i(\alpha_{L^\perp})\leq i(\phi_{L^\perp})\leq \rk L^\perp=1$.

Suppose the lemma holds for $p-1$, and let $L'\subset T_f^\perp M$ be any
subbundle of rank $p-1$ parallel along $D$ with $L\subset L'$. By
\cref{ruled} there is a $D^{n-1}$-ruled isometric immersion
$f':M^n\to\R^{n+p-1}$ whose second fundamental form is
$\sigma'\circ\alpha_{L'}$, where $\sigma':L'\to T_{f'}^\perp M$ is a
parallel bundle isometry.

Now, choosing $L'$ such that
$\Sigma(x)\not\subset L'(x)^\perp$,
by the inductive hypothesis $f'$ is a composition of $\hat f$ near $x$,
i.e., $f'=h'\circ \hat f$ near $x$ for some local isometric immersion
$h'$ between $\R^{n+\ell}$ and $\R^{n+p-1}$. If we further choose $L'$ in
such a way that $\Sigma(x)\not\subset L'(x)$, then $f$ is a
composition of $f'$ near $x$, $f=h''\circ f'$, again by \pref{comp2}.
We conclude that $f= h'' \circ f' = (h''\circ h')\circ \hat f$ is
also a composition of $\hat f$ near $x$.
\qed

\begin{corollary}\po\label{comp}
On each connected component $U$ of an open dense subset of
$M^n$, $f$ is a composition of $f_D$.
\end{corollary}
\proof
Consider an open dense simply connected subset where $\ell_D$
is locally constant, and work on a connected component $V$ of it where
$f_D$ exists by \cref{ruled}. Along the open subset of $V$ where
$\im\alpha\not\subset L_D$, the corollary follows from \pref{ruled2}
applied to $L=L_D$. On the other hand, if $\im \alpha \subset L_D$ along
a connected open subset $U\subset V$, Codazzi equation for
$\eta\in L_D^\perp$ gives
$A_{(\nabla^\perp_Z\eta)_{L_D}}X=A_{(\nabla^\perp_X\eta)_{L_D}}Z=0$.
That is, $L_D^\perp$ is a parallel normal subbundle. Since
$\im \alpha \subset L_D$, we have that $L_D^\perp$ is actually constant
in $\R^{n+p}$. Thus $f(U)\subset\R^{n+\ell_D}$ and the result
also follows on $U$.
\qed

\begin{remark}\po\label{nodirect} {\rm
Observe that in the proof of \pref{ruled2} we did not apply \pref{comp2}
directly to $f$ and $f_D$, but instead inductively. This is so because
all our isometric extensions are extensions by relative nullity: applying
directly \pref{comp2} would give an isometric immersion $h$ with relative
nullity of codimension one only, while the relative nullity of $h$ in
\pref{comp} generically has codimension $p-\ell_D$. The reader should
take this into consideration when trying to apply our results to
submanifolds that are already ruled with big rulings.}
\end{remark}

\begin{corollary}\po\label{extendruled}
If $\hat f$ in \cref{ruled} is a genuine deformation of $f$, then $M^n$
is flat, $L_D=\hat L_D=0$, $\im \alpha \subset L$, and
$D^{n-1}={\cal N}(\alpha)$ almost everywhere. Moreover, $f$ and $\hat f$
singularly extend isometrically along each connected component of an open
dense subset of $M^n$ and, in particular, $\hat f$ is nowhere a strongly
genuine deformation of $f$.
\end{corollary}
\proof
By \cref{comp} we only need to prove the last assertion. Since $\im
\alpha \subset L$, we have that $\phi_\sigma(TM,TM)=0$. Since $f$ and
$\hat f$ are nowhere totally geodesic, by \tref{MainTheorem} we
singularly extend them isometrically using any vector field in $M^n$ not
in $D^{n-1}$.
\qed
\medskip

We point out that all results obtained until now remain valid when the
ambient space is the simply connected space form $\Q_c^m$ of constant
sectional curvature $c$, just by using the exponential map of $\Q^m_c$
when constructing the extensions, e.g., as in (\ref{ExtNat}).

\section{Global applications}\label{ga} 
\hskip 0.6cm The purpose of this section is to give the proof of
\tref{CompactTheorem} and its corollaries. To do this, we use compactness
to transport information along the leaves of relative nullity to the
whole manifold. The use of the intersection of relative nullities makes
the proof short and straightforward, even in the hypersurface case of the
original Sacksteder's theorem, without the need of inductive arguments
or case by case analysis.

\subsection{The intersection of the relative nullities}
\hskip 0.6cm First, we establish some well-known properties of the
splitting tensor adapted to our problem. Let $M^n$ be a Riemannian
manifold and $D$ a smooth totally geodesic distribution on $M^n$. The
{\it splitting tensor} $C$ of $D$ is the map
$C: D \times D^\perp \rightarrow D^\perp$ defined by
$$
C_Y X := C(Y,X) = - (\nabla_X Y)_{D^{\perp}}.
$$
Let $f:M^n \rightarrow \Q^{n+p}_c$ be an isometric immersion of $M^n$
with second fundamental form $\alpha$ and suppose further that $D$ is
contained in the relative nullity ${\cal N}(\alpha)$ of $f$.
Let $\gamma: [0,b] \rightarrow M^n$ be a geodesic such that
$\gamma([0,b))$ is contained in a leaf of $D$. Using the curvature tensor
of $\Q^{n+p}_c$ we easily see that $C_{\gamma'}$ satisfies the Riccati
type ODE
\begin{equation}\label{ambientSplit}
C_{\gamma'}' = C_{\gamma'}^2 + cI,
\end{equation}
where we denote with a ${}'$ the covariant derivative with respect to the
parameter of $\gamma$.

Recall that the shape operator of $f$ in the direction
$\xi \in T_f^{\perp}M$, denoted by $A_{\xi}$, is defined as
$\la A_{\xi}X,Y\ra=\la\alpha(X,Y),\xi\ra$ for $X,Y \in TM$.
For all $Y \in D$ we easily obtain from
Codazzi equation that $\nabla_YA_\xi=A_\xi C_Y+A_{\nabla_Y^\perp\xi}$,
where we understand the operators restricted to $D^\perp$.
If $\xi$ is parallel along $\gamma$ this reduces to
\begin{equation} \label{ODE}
A_{\xi}' = A_{\xi} \circ C_{\gamma '}.
\end{equation}

We will use the splitting tensor of the intersection $\Delta_0$ of the
relative nullities of two isometric immersions, i.e.,
$\Delta_0={\cal N}(\beta)$ for $\beta=\phi_0=\alpha\oplus\hat\alpha$,
which is \eqref{deltatau} for $\tau=0$. We thus need the following
two results for $\Delta_0$.

\begin{lemma}\po \label{totGeoNul}
Let $f:M^n\rightarrow\Q^{n+p}_c$ and $\hat{f}:M^n\rightarrow\Q^{n+q}_c$
be isometric immersions of a Riemannian manifold $M^n$.
Then, along each connected component $U$ of an open dense subset of $M^n$
where $\nu_0=\dim\Delta_0$ is constant, $\Delta_0$ is an integrable
distribution with totally geodesic leaves in $M^n$, $\Q^{n+p}_c$ and
$\Q^{n+q}_c$. In particular, there is a splitting tensor
associated~to~$\Delta_0$ on $U$.
\end{lemma}
\proof
By \pref{mainprops}, $\beta$ is a Codazzi tensor. So
taking $X,Z \in \Delta_0$ and $Y \in TM$ we get
$\beta(Y,\nabla_XZ)=-(\bar{\nabla}^\perp_X \beta)(Y,Z)
=-(\bar\nabla^\perp_Y\beta)(X,Z)=0$, and the lemma follows.
\qed

\begin{lemma}\po \label{extODE}
Let $U \subset M^n$ be an open subset where $\nu_0$ is constant, $\gamma:
[0,b] \rightarrow M^n$ a geodesic with $\gamma([0,b))$ contained in a
leaf of $\Delta_0$ in $U$ joining $x=\gamma(0)$ and $y=\gamma(b)$,
and $P_{\gamma}$ the parallel transport along
$\gamma$ beginning at $t=0$. We have that $\Delta_0(y) =
P_{\gamma}(\Delta_0(x))(b)$, and the splitting tensor
$C_{\gamma'}$ of $\Delta_0$ smoothly extends to $t=b$. In particular, the
ODE (\ref{ODE}) holds up to time $t=b$. Moreover,
$RE(\beta(y))=P_{\gamma}(RE(\beta(x)))(b)$ and
$i(\beta)(y)=i(\beta)(x)$.
\end{lemma}
\proof
Let $J: \Delta_0^{\perp}(\gamma) \rightarrow \Delta_0^{\perp}(\gamma)$ be
the unique solution in $[0,b)$ of the ODE
\begin{equation} \label{(endoJ}
J' + C_{\gamma'} \circ J = 0, \ \ J(0) = I.
\end{equation}
 From (\ref{ambientSplit}) it follows that $J$ also satisfies the linear
ODE with constant coefficients $ J'' + cJ = 0,$ and hence it extends
smoothly to $t=b$, where it is defined in
$P_{\gamma}(\Delta_0^{\perp}(x))(b)$.

For any pair of vector fields $X\!\in\!TM$ and $V\!\in\!\Delta_0^{\perp}$
parallel along $\gamma$, since $\beta$ is Codazzi we have
$$
\bar{\nabla}_{\frac{d}{dt}}^\perp
(\beta(X,J(V)))=\beta(X,(J'+C_{\gamma'}\circ J)(V))=0.
$$
Thus $\beta(X,J(V))$ is parallel along $\gamma$. Since $X(0)$ is
arbitrary, $J$ is invertible in $[0,b]$.
Moreover, since $P_{\gamma}(\Delta_0(x))(b) \subset \Delta_0(y)$ by
continuity, it follows that
$P_{\gamma}(\Delta_0^{\perp}(x))(b)=\Delta_0^{\perp}(y)$.
We conclude that $C_{\gamma'}$ extends smoothly to $[0,b]$ as
$C_{\gamma'} = -J'\circ J^{-1}$ by (\ref{(endoJ}).\linebreak
The last two assertions follow from the parallelism of $\beta(X,J(V))$.
\qed


\subsection{Proofs of the global statements}

\hskip 0.6cm The only ingredient we need to easily obtain
\tref{CompactTheorem} from our local statements is the following.

\begin{proposition}\po\label{prop11}
Let \immr{n}{p} and \immrw{n}{q} be isometric immersions of a compact
Riemannian manifold $M^n$ with $p+q<n$. Then, at each point in $M^n$,
either $i(\beta)\leq p+q-3$, or there are unit vectors
$\xi\in T_f^{\perp}M$ and $\hat\xi\in T_{\hat f}^\perp M$ such that
$\hat A_{\hat\xi}=A_\xi$.
Moreover, the second case holds globally if $\min\{p,q\}\leq 5$.
\end{proposition}

\proof
Let $W\subset M^n$ be the open subset where such unit vectors do
not exist, i.e., where the metric in
$\im(\beta)^\perp\subset T_f^{\perp}M\oplus T_{\hat f}^\perp M$
is definite, and $i(\beta)\geq p+q-2$ if $\min\{p,q\}\geq 6$.
We claim first that $\nu_0>0$ on $W$.

At a fixed a point in $W$, $\beta$ is nondegenerate since
$\im(\beta)^\perp$ is definite. If $\min\{p,q\}\leq 5$, the claim is just
Theorem 3 in \cite{DFComposition}. If otherwise, this is actually the
easiest case in the proof of that theorem for which no hypothesis on
$\min\{p,q\}$ is needed. Indeed, if $X\in RE(\beta)$ we have that
$\dim \im(\beta^X)\cap \im (\beta^X)^\perp
\leq\dim\im(\beta^X)^\perp=p+q-\dim\im(\beta^X)\leq2$.
In this situation, Lemma~6 in \cite{DFComposition}
easily implies that there is $Y\in RE(\beta)$
such that $\Delta_0=\ker(\beta^X)\cap\ker(\beta^Y)=\ker(\lambda)$,
for $\lambda:=\beta^Y|_{\ker(\beta^X)}$. But since $\beta$ is flat
we have that $\im(\lambda)\subset \im(\beta^X)^\perp$, and therefore
$\nu_0=\dim\ker(\lambda)=n-\dim\im(\beta^X)-\dim\im(\lambda)
\geq n-p-q>0$, as claimed.

Let $W'\subset W$ be the open subset where $\nu_0>0$ is minimal in $W$,
and $\gamma\subset W'$ a maximally defined unit geodesic contained in a
maximal leaf of $\Delta_0$ in $W'$. Since, by \lref{totGeoNul}, $\gamma$
is mapped onto a straight line by both $f$ and $\hat f$ and $M^n$ is
compact, $\gamma$~must be defined in a bounded interval $(a,b)$. By
\lref{extODE} the values of $\nu_0$ and $i(\beta)$ are constant along
$\gamma$ up to $t=b$, so $y:=\gamma(b)\not\in W$. Hence, since
$i(\beta)(y)=i(\beta)(\gamma(0))\geq p+q-2$, there are unit vectors
$\xi_0\in T^\perp_{f(y)}M$ and $\hat \xi_0\in T^\perp_{\hat f(y)}M$
such that
$A_{\xi_0}=\hat{A}_{\hat\xi_0}$. If $\xi$ and $\hat\xi$ are their
parallel transports along $\gamma$, by uniqueness of the solutions of
the extended ODE \eqref{ODE} obtained in \lref{extODE} we get
$A_{\xi}=\hat{A}_{\hat\xi}$ also along the whole $\gamma\subset W$, which
contradicts the definition of $W$. We conclude that $W$ is empty.
\qed

\begin{remark}\po\label{rmk03}{\rm
Observe that in the proof above we only used the non-existence of an
unbounded geodesic contained in $\Delta_0$. In particular,
\pref{prop11}, and thus \tref{CompactTheorem} and \cref{DGGeneral}, hold
for complete manifolds if we require that either $f(M)$ or $\hat{f}(M)$
contains no complete straight line instead of compactness.}
\end{remark}

%
%

\noindent {\it Proof of \tref{CompactTheorem}.}
Let $V\subset M^n$ be the open subset where $i(\beta)(x)\geq p+q-2$.
By \pref{prop11} and \cref{par} applied to $\phi_0=\beta$, there exists a
trivially parallel isometry of line bundles parallel along $\Delta_0$,
$\tau:L=\spa\{\xi\}\rightarrow\hat L =\spa\{\hat{\xi}\}$, defined on an
open dense subset $U$ of $V$, and that preserves second
fundamental forms. The result now follows from \cref{CorMain} applied to
each connected component $U'$ of $U$ with this $\tau$, and to
$M^n\setminus \overline V$ with $\tau=0$ since, in either case,
$n+\ell+i(\phi_\tau)\geq n-p-q+3$.
\qed

\medskip

Although \cref{DGGeneral} can be easily proved directly from
\tref{CompactTheorem} and \csref{CorMain} and \ref{par},
we will use the results obtained in \sref{n1}.

\medskip

\noindent {\it Proof of \cref{DGGeneral}.} By \tref{CompactTheorem} we
only need to show that, for $p+q=4$, the immersions singularly extend
isometrically almost everywhere on a subset $U$ where \mbox{$d\geq n-1$}
is constant. Clearly, this is the case if $d=n$, since both immersions
would be totally geodesic in $U$ and we isometrically extend them with
$N^{n+1}=U\times\R$, $F=f\times Id$ and $\hat F=\hat f\times Id$.

If $d=n-1$ on $U$, we have as in the proof of \tref{CompactTheorem} and
again by \cref{CorMain} that $p+q=4$, $\ell=1$, $i(\phi_\tau)=2$, and
$D^{n-1}={\cal N}(\alpha_{L^\perp})\cap{\cal N}(\hat\alpha_{\hat L^\perp})$
along an open dense subset of $U$. If $U'\subset U$ is the open subset
where $L_D$ in \eqref{ld} is nonzero, then $L=L_D$, $\hat L = \hat L_D$,
$f_D=\hat f_D$, and thus by \cref{comp} $f$ and $\hat f$ (regularly)
extend isometrically almost everywhere on~$U'$. On $U\setminus U'$,
$L_D=0$ and $D^{n-1}=\Delta_0$ almost everywhere, so \csref{par} and
\ref{ruled} tell us
that there is an \isim $f':U\to\R^{n+1}$ with second fundamental form
$\alpha_L=\hat\alpha_{\hat L}$. By \cref{extendruled}, the pairs
$\{f,f'\}$ and $\{\hat f,f'\}$ both singularly extend isometrically, and
since the codimension of $f'$ is one, the pair $\{f,\hat f\}$ also
singularly extends isometrically almost everywhere on~$U\setminus U'$.
\qed

%
%
%

\begin{remark}\po{\rm
\cref{DGGeneral} for $p=q=2$ reduces to the main result in
\cite{DGCompact}, except for the fact that singular flat extensions can
occur in the former. This is a consequence of a gap in \cite{DGCompact},
whose long and involved case by case proof did not cover all
possibilities.
}
\end{remark}

\noindent {\it Proof of \cref{theoExtra}.} By Proposition 26 in
\cite{DFGenuines} and \tref{CompactTheorem}, if an isometric immersion
\immrw{n}{q} is a strongly genuine deformation of \immr{n}{p}, then the
$k$-th Pontrjagin form $p_k$ of $M^n$ vanishes for any $k$ such that $4k>
3(p+q-3)$.
\qed

\section{The space forms case}\label{sphereCase} 
\hskip 0.6cm As we pointed out, \tref{CompactTheorem} and Corollary
\ref{DGGeneral} hold for compact manifolds when the ambient space is the
hyperbolic space following the same proofs. In this section we show that
they also hold for complete submanifolds in the sphere under a
mild codimension condition.

\medskip

For the following, recall that $\rho(m)-1$ is the maximum number
of pointwise linearly-independent vector fields on $\mathbb{S}^{m-1}$.

\begin{lemma}\po\label{sphereExt}
Let $f: M^n \rightarrow \mathbb{S}^{n+p}$ be an isometric immersion and
$D^d$ a nontrivial totally geodesic distribution contained in the
relative nullity of $f$. If there exists a nonconstant geodesic $\sigma:
[0, \infty) \rightarrow M^n$ in $D^d$, then the splitting tensor
$C_{\sigma'}$ associated to $D^d$ has no real eigenvalues. In particular,
such a geodesic cannot exist if $\rho(n-d) < d + 1$.
\end{lemma}
\proof
By (\ref{ambientSplit}), $C_{\sigma'}$ is given by
$$
(P_{\sigma}^{-1}\circ C_{\sigma'} \circ P_{\sigma})(t)
=(\sin(t)I+\cos(t)C_{\sigma'}(0))(\cos(t)I-\sin(t)C_{\sigma'}(0))^{-1},
$$
where $P_{\sigma}$ is the parallel transport along $\sigma$. Since
$C_{\sigma'}$ is defined for all $t \geq 0$, we easily conclude that
$C_{\sigma'(0)}$ has no real eigenvalues.

For the last assertion, choose a basis $\{T_1,...,T_d\}$ of $D(x)$.
By the first assertion, for any unit vector $Z \in D^\perp(x)$ and
$a,a_1,...,a_{d} \in \R$, the equation
$0=aZ+\sum_{i=1}^{d}a_iC_{T_i}Z=aZ+C_T Z$
implies that $a=a_i=0$, where $T=\sum_{i=1}^da_iT_i$. Hence
$Z,C_{T_1}Z,...,C_{T_d}Z$ are linearly independent in $D^{\perp}(x)$.
Since this holds for any unit vector $Z \in D^{\perp}(x)$, considering
$Z$ as the position vector of the unit sphere
$\mathbb{S}^{n-d-1}\subset D^{\perp}(x)$ we get $d$ nonvanishing linearly
independent vector fields in $\mathbb{S}^{n-d-1}$. Hence,
$d\leq\rho(n-d)-1$.
\qed
\medskip

\noindent {\it Proof of \tref{GloRes}.} By \lref{sphereExt},
geodesics in $\Delta_0$ cannot be defined for arbitrary large time if
$\mu_n < n - p - q$ when the ambient space is the sphere. Thus, as
observed in Remark \ref{rmk03}, \pref{prop11} holds for complete
manifolds when the ambient spaces are spheres as long as $p+q<n-\mu_n$.
\qed
\medskip

\noindent {\it Proof of \cref{corEsf}.} It is analogous to the
one for \cref{DGGeneral} using \tref{GloRes} instead of
\tref{CompactTheorem}, just observing that for small codimensions
we can simplify the assumptions on $\mu_n$.
\qed
\medskip


\addcontentsline{toc}{chapter}{Bibliography}
\bibliographystyle{abbrv}

{\renewcommand{\baselinestretch}{1}
\hspace*{-20ex}\begin{tabbing}
\indent \= IMPA -- Estrada Dona Castorina, 110 \\
\> 22460-320 --- Rio de Janeiro --- Brazil \\
\> luis@impa.br\ \ \ \ \ -- \ \ \ \ felippe@impa.br \end{tabbing}
}
\end{document}